\documentclass[a4paper, 11pt]{article}
\usepackage{amsmath}
\usepackage{amssymb}
\newcommand{\comment}[1]{}

\newcommand{\mybox}{\hfill \mbox{$\Box$}}

\def \R{{\mathbb R}} 
\def \N{{\mathbb N}}  
\def\square{\hfill \vrule height 1.5ex width 0.8ex depth -.1ex}

\def\taus{\tau_\sigma}
\def\Rc{\overline{\R_{+}}}

\def\K{\mathcal{K}}
\def\S{\mathcal{S}}

\def\Sa{\mathcal{S}_{a}}

\def\T{\mathcal{T}}
\def\Ta{\mathcal{T}_a}

\def\P{\mathcal{P}}
\newtheorem{teo}{Theorem}[section]
\newtheorem{itlemma}{Lemma}[section]
\newtheorem{lema}{Lemma}[section]

\newtheorem{defin}{Definition}[section]
\newtheorem{prop}{Proposition}[section]
\newtheorem{itremark}[itlemma]{Remark}
\newtheorem{as}{Assumption}
\newenvironment{remark}{\noindent\begin{itremark}\rm}{\end{itremark}}

\title{\LARGE \bf
An Extension of LaSalle's Invariance Principle for Switched Systems
}

\comment{
\author{J. L. Mancilla-Aguilar and R.A. Garc\'{\i}a 
\thanks{J.L. Mancilla-Aguilar is with the Department of Mathematics,
    Faculty of Engineering, University of Buenos Aires.
        {\tt\small jmancil@fi.uba.ar}}%
\thanks{R. A. Garc\'{\i}a is with the Department of Physics and
        Mathematics of the Instituto Tecnol\'ogico de Buenos Aires, and
        with the Faculty of Engineering of the University of Buenos
        Aires.
        {\tt\small ragarcia@itba.edu.ar}} }

}

\author{J. L. Mancilla-Aguilar \thanks{ \jma} 
 \;and R.A. Garc\'{\i}a \thanks{\rg}}

\def\jma{Department of Mathematics,
    Faculty of Engineering, University of Buenos Aires.
        {\tt\small jmancil@fi.uba.ar}\ Work partially supported by UBA I039.}

\def\rg{Department of Physics \&
        Mathematics of the Instituto Tecnol\'ogico de Buenos Aires; and
        the Faculty of Engineering of the University of Buenos
        Aires.
        {\tt\small ragarcia@itba.edu.ar}\ Work partially supported by UBA I039.}

\begin{document}
\maketitle

\begin{abstract}
This paper addresses invariance principles for a certain class of switched nonlinear  
systems. We provide an extension of LaSalle's Invariance Principle for these systems and 
state asymptotic stability  criteria. We also present some related results that deal 
with the compactness of the trajectories of these switched systems and that are 
interesting by their own.

\end{abstract}

\section{Introduction}
In recent years, switched systems have deserved a great deal of  attention from the 
Systems Engineering and Computer Sciences communities. In particular, the stability 
properties of the common equilibrium solutions have been intensively investigated, see 
e.g. \cite{DeCarlo}, \cite{Liberzonbook} and \cite{liberzon}, respectively, and the 
references therein. Although switched systems whose component subsystems are autonomous 
(the class of switched systems that we will consider here) are in essence nonautonumous 
systems, and as such were investigated by different authors, (see \cite{Kloeden} and 
references therein), their stability properties can also be studied by means of multiple 
Lyapunov functions, \cite{Bacciotti}, \cite{Branicky}, \cite{Hespanha}, 
\cite{Hespanha2}, \cite{Peleties} (this approach is very attractive since it enables us 
to study their uniform, in the sense of the switching signals considered, stability 
properties). In this context, several Lyapunov-like  results and different invariance 
principles have been recently proposed. Among the invariance  results,  in 
\cite{Hespanha} LaSalle's invariance principle is extended to switched linear systems 
under rather general switching; a LaSalle-like invariance principle for switched 
nonlinear systems under more restrictive switching than that of \cite{Hespanha} is 
proposed in \cite{Bacciotti} and by  using the small-time norm-observability hypothesis, 
theorems for switched nonlinear systems in the same line as those in \cite{Hespanha} are 
proved in \cite{Hespanha2}. A version of LaSalle's invariance principle for 
deterministic hybrid automata with a finite number of discrete states is presented in 
\cite{lygeros}. 

In this paper, we extend LaSalle's invariance principle to switched nonlinear systems 
under mild restrictions in the class of switchings, since we consider switching signals 
with a positive average dwell-time (see the definition below). This extended principle 
enables us to obtain (uniform) asymptotic stability criteria for this class of systems.

The results that we present here enable us not only to extend partially some of those in 
\cite{Hespanha2} and improve those in \cite{Bacciotti}, but get a better comprehension 
of the structure of invariant sets and of the compactness properties of trajectories of 
this class of switched systems. 

The paper is organized as follows. In section 2 we establish the notation, and present 
the basic definitions and main results of the paper. In section 3 we exhibit some 
examples that show the application of our results to systems  to which those mentioned 
above cannot be applied or are of little help. We study the invariance for switched 
systems and prove our Invariance Principle in section 4. Section 5 is devoted to the 
proof of one of the main results. In section 6 we present the conclusions and finally in 
the Appendix we present some results about compactness of the trajectories of the 
switched systems under study. These results, that are used in some proofs along the 
paper, are also important by their own.

\section{Basic definitions and main results}
Throughout, $\R$, $\R_{+}$, $\N$ and $\N_0$ denote the sets of real, nonnegative real, 
natural and nonnegative integer numbers, respectively. 
We use $|\cdot|$ to denote the Euclidean norm on $\R^n$. As usual, by a 
$\mathcal{K}$-function we mean a function
$\alpha:\R_+ \rightarrow \R_+$ that is strictly increasing and continuous,
and satisfies $\alpha(0)=0$, by a $\mathcal{K}_{\infty}$-function one that is in 
addition unbounded, and we let $\mathcal{KL}$ be the class of functions
$\R_+ \times \R_+ \rightarrow \R_+$ which are of class $\mathcal{K}_{\infty}$ on
the first argument and decrease to zero on the second argument. Let 
$\mathcal{C}(\R^n,\R^n)$ denote the set of all the  continuous maps from $\R^n$ to 
$\R^n$.

Given a family $\P=\{f_\gamma\in \mathcal{C}(\R^n,\R^n): \gamma \in \Gamma \}$, where 
$\Gamma=\{1,\ldots,m\}$, we consider the switched system described by 
\begin{eqnarray}
\dot{x}=f(x,\sigma) \label{ss}
\end{eqnarray}
where $x$ takes values in $\R^n$, $\sigma: \R_{+} \rightarrow \Gamma$ is a {\em 
switching signal}, {\em i.e.},  $\sigma$ is piecewise constant and continuous from the 
right and $f:\R^n \times \Gamma \rightarrow \R^n$ is defined by 
$f(\xi,\gamma)=f_\gamma(\xi)$. In what follows we consider that the set $\Gamma$ is 
equipped with the discrete metric; in consequence  $\Gamma$ is a compact metric space 
and $f$ a continuous function.

We will denote by $\S$ the set of all the switching signals. We recall that a piecewise 
smooth curve $x:\mathcal{I} \rightarrow \R^n$,  with $\mathcal{I}=[0,T)$ or $[0,T]$ with 
$0<T \le +\infty$ is a {\em solution} of (\ref{ss}) corresponding to $\sigma \in \S$ if 
$\dot{x}(t)=f(x(t),\sigma(t))$ for all $t \in (t_i,t_{i+1})\cap \mathcal{I}$, where 
$t_0=0$ and $t_1, t_2,\ldots$ are the consecutive discontinuities (switching times) of 
$\sigma$. 

A pair $(x,\sigma)$ is a {\em trajectory} of (\ref{ss}) if $\sigma \in \S$ and $x$ is a 
solution of (\ref{ss}) corresponding to $\sigma$. We say that a trajectory $(x,\sigma)$ 
of (\ref{ss}) is {\em maximal} and denote its domain by 
$\mathcal{I}_{(x,\sigma)}=[0,t_{(x,\sigma)})$, if $x$ is a maximal solution of 
(\ref{ss}) corresponding to $\sigma$. We observe that, due to standard results on 
ordinary differential equations, either $t_{(x,\sigma)}=\infty$ or 
$t_{(x,\sigma)}<\infty$ and $x$ is unbounded. Let $\T$ denote the set of all the maximal 
trajectories of (\ref{ss}). 

Since in many applications, the admissible switching signals, and consequently the 
admissible trajectories of the switched system, are not completely arbitrary since they 
are subjected to constraints which may concern their functional nature or be dependent 
on the states $x(t)$ of the system (see \cite{Liberzonbook, liberzon}), we will suppose 
that the class of admissible trajectories of (\ref{ss}) is a sub-family $\T^\prime$ of 
the whole family of trajectories $\T$.

In order to take into account some kind of state-dependent constraints on the admissible 
trajectories, as done in  \cite{Hespanha} we introduce the following family of 
trajectories: Given a covering $\chi:=\{\chi_\gamma:\;\gamma \in \Gamma\}$ of $\R^n$, 
i.e., $\chi_\gamma \subseteq \R^n$ and $\R^n=\bigcup_{\gamma \in \Gamma}\chi_\gamma$, we 
define $\T[\chi]$ as the set of all the trajectories $(x,\sigma) \in \T$ which verify 
the condition: $x(t)\in \chi_{\sigma(t)}$ for all $t \in \mathcal{I}_{(x,\sigma)}$. We 
say that $\chi$ is a closed covering if $\chi_\gamma$ is closed for every $\gamma \in 
\Gamma$.

In what follows we assume that the following standing hypothesis holds
\begin{as} There exists a closed covering $\chi$ of $\R^n$ such that $\T^\prime 
\subseteq \T[\chi]$. \mybox
\end{as}

Since we are mainly interested in the stability analysis of the zero solutions of system 
(\ref{ss}), we will also suppose the following hypothesis holds 
\begin{as} $f(0,\gamma)=0$ for all $\gamma \in \Gamma^*$, where  $\Gamma^*=\{\gamma \in 
\Gamma:\;0 \in \chi_\gamma\}$;\mybox
\end{as}
and adopt the following definitions of stability. 
\begin{defin} \label{stability} We say that a family $\T^\prime$ of maximal trajectories 
of (\ref{ss}) is 
\begin{enumerate}
\item
{\em uniformly stable} if there exists a function $\alpha \in \K_{\infty}$ such that 
for every $(x,\sigma) \in \T^\prime$, 
\begin{eqnarray*} 
|x(t)| \le \alpha(|x(t_0)|)\quad \forall t\ge t_0,\;\forall t_0\ge 0.
\end{eqnarray*}     
\item
{\em Globally asymptotically stable} if it is uniformly stable and, in addition, for 
every trajectory $(x,\sigma)\in \T^\prime$, $x(t)$ converges to $0$ as $t \rightarrow 
\infty$.
\item
{\em Globally uniformly asymptotically stable} if there exists a function $\beta \in 
\mathcal{KL}$ such that for every $(x,\sigma) \in \T^\prime$,
 \begin{eqnarray*} 
|x(t)| \le \beta(|x(t_0)|,t-t_0)\quad \forall t\ge 0, \;\forall t_0\ge 0.
\end{eqnarray*} 
\end{enumerate} \mybox
\end{defin}
\begin{remark}\label{gas}
With the same technique used to prove Proposition 2.5 of \cite{Lin}, it can be
shown that the definition of global uniform asymptotic stability given above is 
equivalent to the following (more classical) one:\\
A family $\T^\prime$ of maximal trajectories of (\ref{ss}) is globally uniformly 
asymptotically stable if it is uniformly stable and 
\begin{itemize}
\item[(a)]
for each $R>0$ and each $\varepsilon > 0$, there exists
$T\ge 0$ such that for all $(x, \sigma) \in \T^\prime$ and all $t_0 \ge 0$
\[ |x(t_0)| < R \Longrightarrow |x(t)|< \varepsilon \quad \forall t \ge t_0+T.\]
\end{itemize}
\end{remark} \mybox

Several Lyapunov-like theorems which involve the use of multiple Lyapunov functions 
(see, among others, \cite{Branicky}, \cite{DeCarlo}) allow us to establish the stability 
or asymptotic stability of a family of admissible trajectories. The following one which 
is based on results given in \cite{Branicky}, \cite{Peleties} (see also \cite{Hespanha}) 
is an example of such theorems. It is convenient to introduce here the following
\begin{defin} \label{wLl} A function $V: \R^n \times \Gamma \rightarrow \R$ is a weak 
Lyapunov-like function for the family $\T^\prime$ if it is continuously differentiable 
with respect to the first argument and verifies
\begin{enumerate}
\item there exist
$\alpha_1$ and $\alpha_2$ of class $\K_\infty$ so that $\alpha_1(|\xi|)\le 
V(\xi,\gamma)\le \alpha_2(|\xi|)$ for all $(\xi, \gamma) \in \R^n \times \Gamma$ such 
that $\xi \in \chi_\gamma$;
\item $\frac{\partial V}{\partial \xi}(\xi,\gamma)f(\xi,\gamma) \le 0$, for all 
$(\xi,\gamma) \in \R^n \times \Gamma$ such that $\xi \in \chi_{\gamma}$;
\item for every trajectory $(x,\sigma)\in \T^\prime$ and any pair $t_i < t_j$ of 
switching times such that $\sigma(t_i)=\sigma(t_j)$, $V(x(t_{j}),\sigma(t_{j}))\le 
V(x(t_{i+1}),\sigma(t_i))$.
\end{enumerate} \mybox
\end{defin}
\begin{teo} \label{strict} \rm  Suppose there exists a weak Lyapunov-like function $V$ 
for $\T^\prime$. Then $\T^\prime$ is uniformly stable. 

If, in addition, $V$ verifies
\begin{itemize}
\item[2$^\prime$.] there exists a positive-definite function $\alpha_3$ such that
$\frac{\partial V}{\partial \xi}(\xi,\gamma)f(\xi,\gamma) \le -\alpha_3(|\xi|)$, for all 
$(\xi,\gamma) \in \R^n \times \Gamma$ such that $\xi \in \chi_{\gamma}$,
\end{itemize}
then $\T^\prime$ is globally uniformly asymptotically stable. \mybox
\end{teo} 

As was pointed out above, this work is concerned with invariance principles for switched 
systems and, in particular, with extensions of LaSalle's invariance principle to this 
class of systems. In this regard we will show, under suitable hypotheses, that the 
existence of a weak Lyapunov-like function $V$ for $\T^\prime$, allows us to obtain 
conclusions about the asymptotic behavior of a bounded solution $x$ of (\ref{ss}) 
corresponding to some switching signal $\sigma$ so that $(x,\sigma)\in \T^\prime$, and, 
further, to obtain some asymptotic stability criteria. 

As was discussed in \cite{Hespanha} and also in \cite{Bacciotti}, in order to obtain 
LaSalle-like asymptotic stability criteria by exploiting the knowledge of a weak 
Lyapunov-like function $V$, some form of regularity in the switching signals regarding 
the distance between consecutive switching times is needed. In this paper we will 
consider switching signals which have a positive average dwell-time, more precisely, 
\begin{defin} We say that the switching signal $\sigma$ has an average dwell-time 
$\tau_D>0$ and a chatter bound $N_0 \in \N$ if the number of switching times of $\sigma$ 
in any open finite interval $(\tau_1,\tau_2) \subset \R_{+}$ is bounded by 
$N_0+(\tau_2-\tau_1)/\tau_D$. \mybox
\end{defin}     

We denote by $\Sa[\tau_D,N_0]$ the set of all the switching signals which have an 
average dwell-time $\tau_D>0$ and a chatter bound $N_0 \in \N$ and by $\Ta[\tau_D, N_0]$ 
the subclass of all the trajectories of (\ref{ss}) corresponding to some $\sigma \in 
\Sa[\tau_D,N_0]$. Let $\Sa=\bigcup_{\tau_D >0, N_0 > 0} \Sa[\tau_D,N_0]$ and let $\Ta$ 
denote the subclass of all the trajectories of (\ref{ss}) corresponding to some $\sigma 
\in \Sa$, i.e. $\Ta=\bigcup_{\tau_D>0, N_0 > 0} \Ta[\tau_D, N_0]$.

We note that the set of switching signals $\sigma$ which have a dwell-time $\tau_D>0$, 
i.e., $\inf_{k\ge 0}t_{k+1}-t_k \ge \tau_D$, is a subset of $\S[\tau_D,1]$.

From now on we suppose that the following additional hypothesis holds.
\begin{as} $\T^\prime \subseteq \Ta$. \mybox
\end{as}

In order to establish the main results of this paper, we need to introduce 
some more trajectory families.

Given a continuous function $V:\Omega \times \Gamma \rightarrow \R$, with $\Omega$ an 
open subset of $\R^n$, we consider the following families of trajectories associated 
with $V$.

$\T_V$ is the class of all  the trajectories $(x,\sigma) \in \T$ which verify the 
conditions:
\begin{enumerate}
\item $x(t) \in \Omega$ for all $t \in \mathcal{I}_{(x,\sigma)}$;
\item for any pair of times $t, t^\prime \in \mathcal{I}_{(x,\sigma)}$ such that $t \le 
t^\prime$ and $\sigma(t)=\sigma(t^\prime)$, $V(x(t),\sigma(t))\ge 
V(x(t^\prime),\sigma(t^\prime))$. 
\end{enumerate}

$\T_V^{*}$ is the sub-family of $\T_V$ whose members $(x,\sigma)$ verify the condition
\[ \sigma(t)=\sigma(t^\prime) \Longrightarrow V(x(t),\sigma(t))= 
V(x(t^\prime),\sigma(t^\prime)).\] 
\begin{remark} If $V$ is a weak Lyapunov-like function for $\T^\prime$ it readily 
follows that $\T^\prime \subseteq \T_V$. \mybox
\end{remark} 

Finally, we introduce the following notion of weak-invariance for nonempty subsets of 
$\R^n \times \Gamma$.
\begin{defin} Given a family $\T^*$ of maximal trajectories of (\ref{ss}), we say that a 
nonempty subset $M \subseteq \R^n \times \Gamma$ is weakly-invariant with respect to 
$\T^*$ if for each $(\xi,\gamma)\in M$ there is a trajectory $(x,\sigma) \in \T^*$ such 
that $x(0)=\xi$, $\sigma(0)=\gamma$ and $(x(t),\sigma(t)) \in M$ for all $t \in 
\mathcal{I}_{(x,\sigma)}$. \mybox
\end{defin}

Now we are in position to state the following asymptotic stability criterion, which is 
one of our main results.
\begin{teo} \label{cordin} \rm  Suppose that there exists a weak Lyapunov-like function 
$V$ for $\T^\prime$ such that $M=\{0\}\times \Gamma^*$ is the maximal weakly-invariant 
set w.r.t. $\T_V^*\cap \T[\chi]\cap\Ta$. 

Then $\T^\prime$ is globally asymptotically stable. If, in addition, $\T^\prime 
\subseteq \Ta[\tau_D,N_0]$ for some $\tau_D>0$ and some $N_0\in \N$ then $\T^\prime$ is 
globally uniformly asymptotically stable.  \mybox
\end{teo}
\begin{remark} \rm From the proof of Theorem \ref{cordin}, which can be found in section 
5, it follows that in the case when $\T^\prime \subseteq \Ta[\tau_D,N_0]$ for some 
$\tau_D>0$ and some $N_0\in \N$, the thesis of the theorem still holds if one assumes 
the weaker hypothesis $M=\{0\}\times \Gamma^*$ is the maximal weakly-invariant set 
w.r.t. $\T_V^*\cap \T[\chi]\cap\Ta[\tau_D,N_0]$. \mybox
\end{remark}   
The following result can be readily deduced from Theorem \ref{cordin}.
\begin{teo} \label{sontagteo} \rm  Suppose that there exists a weak Lyapunov-like 
function $V$ for $\T^\prime$. Suppose, in addition, that there exists a family 
$\{W_\gamma : \R^n \rightarrow \R, \gamma \in \Gamma\}$ of continuous and nonnegative 
definite functions such that
\begin{enumerate}
\item $\frac{\partial V}{\partial \xi}(\xi,\gamma)f(\xi,\gamma) \le -W_\gamma(\xi)$, for 
all $(\xi,\gamma) \in \R^n \times \Gamma$ such that $\xi \in \chi_\gamma$;
\item for each $\gamma \in \Gamma$, the system  
\begin{eqnarray}
\label{dis}
\dot{x} &=& f_\gamma(x), \quad y = W_\gamma(x), 
\end{eqnarray}
is {\em zero small-time distinguishable}. (We recall that a systems (\ref{dis}) is zero 
small-time distinguishable,  if for every  $\delta>0$, $x(0) = 0$ whenever 
$W_{\gamma}(x(t))=0$ for all $t \in [0, \delta)$).
\end{enumerate}
Then $\T^\prime$ is globally asymptotically stable. If, in addition, $\T^\prime 
\subseteq \Ta[\tau_D,N_0]$ for some $\tau_D>0$ and some $N_0\in \N$ then $\T^\prime$ is 
globally uniformly asymptotically stable. \mybox
\end{teo}

{\em Proof.} Let $M$ be the maximal weakly-invariant set w.r.t. $\T_V^*\cap 
\T[\chi]\cap\Ta$.  In order to prove the theorem, it suffices to show that $M = \{0\} 
\times \Gamma^*$, since then the hypotheses of Theorem \ref{cordin} will be fulfilled. 
Note first that $\{0\} \times \Gamma^*$ is weakly-invariant w.r.t.$\T_V^*\cap 
\T[\chi]\cap\Ta$ due to Assumption 2 and the fact that, due to 1. in Definition 
\ref{wLl}, $V(0, \gamma) = 0 \, \forall\, \gamma \in \Gamma^*$.

Let then $(\xi, \gamma) \in M$; it follows that there exists a trajectory $(x,\sigma) 
\in \T_V^{*}\cap\T[\chi]\cap \Ta$ such that $(x(0),\sigma(0))=(\xi,\gamma)$ and 
$(x(t),\sigma(t)) \in M$ for all $t \ge 0$. Let $\tau > 0$ such that $\sigma(t) = 
\gamma$ for all $t \in [0, \tau)$. It follows from the definition of $\T_V^{*}$ that 
$\frac{\partial V}{\partial \xi}(x(t), \sigma(t)) f(x(t), \sigma(t))= 0$ for all $t \in 
[0, \tau)$ and hence $W_{\gamma}(x(t)) = 0$ for every $t \in [0, \tau)$. From the zero 
small-time distinguishability assumption, we have that $\xi=x(0) = 0$. Since $(x,\sigma) 
\in \T[\chi]$, $x(0) \in \chi_{\sigma(0)}$ and consequently $\gamma \in \Gamma^*$.\\
\makebox{} \square
\begin{remark}\rm Theorem \ref{sontagteo} is a partial generalization of Theorem 7 in 
\cite{Hespanha2}, since the zero small-time distinguishability hypothesis is weaker than 
the small-time norm-observability Assumption 2. in \cite{Hespanha2} and since we obtain 
{\em uniform} asymptotic stability in the case when $\T^\prime \subseteq 
\Ta[\tau_D,N_0]$ for some $\tau_D>0$ and some $N_0\in \N_0$, but our hypothesis 
$\T^\prime \subseteq \Ta$ is slightly stronger than the hypothesis about the regularity 
of the switching signals considered in \cite{Hespanha2} (see Assumption 3. of that 
paper). \mybox
\end{remark}
The proof of Theorem \ref{cordin} is based on the following extension of the well known 
invariance principle for dynamical systems described by differential equations of 
LaSalle (see \cite{lasalle}) to switched systems. This extension  is other of the main 
results of this work. 

Let $\pi_1:\R^n \times \Gamma \rightarrow \R^n$ be the projection onto the first 
component.
\begin{teo} \label{Lasalle}
\rm Let $\chi$ be a closed covering of $\R^n$ and let $V:\Omega \times \Gamma 
\rightarrow \R$, with $\Omega$ an open subset of $\R^n$, be continuous. Suppose that 
$(x,\sigma)$ is a trajectory belonging to $\T_V \cap \T[\chi] \cap \Ta[\tau_D,N_0]$ such 
that for some compact subset $B \subset \Omega$, $x(t) \in B$ for all $t \ge 0$. Let $M 
\subseteq \R^n \times \Gamma$ be the largest weakly-invariant set w.r.t. $\T_V^{*}\cap 
\T[\chi]\cap \Ta[\tau_D,N_0]$ contained in $\Omega \times \Gamma$. 

Then, $x(t)$ converges to $\pi_1(M)$ as $t \rightarrow \infty$. \mybox
\end{teo}
\begin{remark}\label{bacc}  In the case when we restrict the hypotheses of Theorem 
\ref{Lasalle} to those of Theorems 1 and 2 in \cite{Bacciotti}, we obtain more precise 
results related to the size of the attracting sets involved.  We shall prove our 
assertion for Theorem 1 in \cite{Bacciotti} only, since the proof for the other is 
similar. In what 
follows, and in order to prove our claim, we refer to the notation and definitions of 
that paper. Let $x(t)$ be a dwell-time solution (in the sense of \cite{Bacciotti}) with 
initial condition $x(0)\in \Omega_l$, a dwell-time $\tau_D > 0$ and 
generated by a switching signal $\sigma$ in a switched system that admits a common weak 
Lyapunov function (in the sense above) $V:\Omega \rightarrow \R_{+}$. Let 
$W:\Omega_l\rightarrow \R_+$ be the restriction to $\Omega_l$ of the function $V$.  It 
is easy to see that the trajectory $(x, \sigma)$ belongs to $\T_W \cap \Ta[\tau_D, 1]$, 
and that there exists a compact set $B\subset \Omega_l$ so that  $x(t) \in B$ for all $t 
\ge 0$. Therefore the trajectory  $(x,\sigma)$ verifies the hypotheses of Theorem 
\ref{Lasalle} (with $W$ in place of $V$ and the trivial covering of $\R^n$, 
$\chi_\gamma=\R^n$ for all $\gamma \in \Gamma$). Hence $x(t)$ converges to $\pi_1(M)$ as 
$t \rightarrow \infty$, where $M$ is the largest weakly-invariant set w.r.t. 
$\T_W^*\cap\Ta$ contained in $\Omega_l \times \Gamma$.

On the other hand, Theorem 1 asserts that $x(t)$ converges to $M^\prime$, where 
$M^\prime$ is the union of all the compact, weakly-invariant sets (in the sense of 
\cite{Bacciotti}) which are contained in $Z \cap \Omega_l$. We will prove that $\pi_1(M) 
\subseteq M^\prime$ and therefore our assertion about the sizes of the attracting sets.

Pick $\xi \in \pi_1(M)$. Then there exist $\gamma \in \Gamma$ and a trajectory $(x^*, 
\sigma^*) \in  \T_W^* \cap \Ta$ 
such that $(x^*(0), \sigma^*(0)) = (\xi, \gamma)$ and $(x^*(t), \sigma^*(t)) \in M$ for 
all $t \in \mathcal{I}_{(x^*,\sigma^*)}$. In consequence $\sigma^*(t) = \gamma$ for all 
$t \in [0, \delta]$ with $0 < \delta < t_1$ and $t_1$ the first switching time of 
$\sigma^*$. It follows from the definition of $\T_W^*$ that $\frac{\partial
V}{\partial \xi}(x^*(t)) f(x^*(t), \gamma) = 0$ and hence that $x^*(t) \in Z \cap 
\Omega_l$ for all $t\in [0,\delta]$. Next, $x^*([0,\delta])$ is a compact 
weakly-invariant set contained in $Z\cap\Omega_l$ and consequently $\xi \in M^\prime$. 

It must be remarked that the attracting set $\pi_1(M)$ corresponding to the 
weakly-invariant set considered in Theorem \ref{Lasalle}, may be considerably smaller 
than the attracting sets given in Theorems 1 and 2 of \cite{Bacciotti}, as we exhibit in 
Example 2 below. \\
\makebox{}\mybox
\end{remark}   
\section{Examples}
{\bf Example 1.}  Consider the switched system in $\R^2$ given by the family 
$\{f_1,f_2\}$, with  
\begin{eqnarray*}
f_1(\xi) = \left(\begin{array}{c}-2\xi_1 - 2\xi_2 \\ 2 \xi_1 \end{array} \right),\;\; 
f_2(\xi) = \left(\begin{array}{c}-\xi_2\\ \xi_1 \end{array} \right). 
\end{eqnarray*}
Let $\T^\prime$ be the set of all the maximal trajectories $(x,\sigma)$ whose switching 
signals $\sigma$ are given by the feedback rule  
\begin{eqnarray*}
\sigma(t)= \left \lbrace \begin{array}{rrr} 1 &\makebox{if} & x_1(t) < 0 \\ 2 
&\makebox{if} &x_1(t) \geq 0. \end{array} \right.
\end{eqnarray*}
Observe that the origin  is a stable focus for the first subsystem and a center for the 
other, and that the trajectories of both are running counterclockwise.
 
Since the time needed by any nontrivial trajectory of the subsystem $\dot{x}=f_1(x)$ 
($\dot{x}=f_2(x)$) to go from the positive (resp. negative) $x_2$ axis to the negative 
(resp. positive) one is constant, clearly $\T^\prime \subseteq \Ta[\tau_D, 1]$ for some 
$\tau_D > 0$. If we consider the closed covering of $\R^2:\; \chi_1 = \{\xi : \xi_1 \leq 
0\}, \, \chi_2 = \{\xi : \xi_1 \geq 0\}$, then $\T^\prime \subseteq \T[\chi] \cap 
\Ta[\tau_D, 1]$.

The function $V : \R^2 \times \{1, 2\} \rightarrow \R$ defined by $V(\xi, i) = |\xi|^2$ 
is clearly a weak Lyapunov-like function for $\T^\prime$. We claim that $M = \{0\} 
\times\{1, 2\}$. In fact, let $(\xi,\gamma) \in M$. Then there exists $(x, \sigma) \in 
\T_V^*\cap\T[\chi]\cap \Ta$ such that $x(0)=\xi$, $\sigma(0)= \gamma$ and 
$(x(t),\sigma(t))\in M$ for all $t\ge 0$. From the facts that $x(t)$ cannot remain 
forever in the right half-plane, where $V(x(t), \sigma(t))$ is constant, that $V(x(t), 
\sigma(t))$ is strictly decreasing when $x(t)$ is in the open left-half plane, and from 
the 
definition of $\T_V^{*}$, it follows readily that $(x(t), \sigma(t))$ cannot belong to 
$M$ unless $x(t) = 0$ for all $t \geq 0$, and the claim follows. Hence, according to 
Theorem \ref{cordin}, $\T^\prime$ is globally uniformly asymptotically stable.\\
\\
{\bf Example 2.} Consider now the two systems in $\R^2$ given by:
\begin{align*}
\dot{x} = f_1(x) = \left(\begin{array}{c}-x_1 - x_2 \\ x_1 \end{array} \right) 
\,\makebox{and}\;\, \dot{x} =  f_2(x) = -\frac{x}{1 + |x|^4}.  
\end{align*}
Let $W_1(\xi) = \xi_1^2$, $ W_2(\xi) = \frac{|\xi|^2}{1 + |\xi|^4}$ 
and $V(\xi, i)  = |\xi|^2/2, \;i = 1, 2$. Then for every $\xi \in \R^2$, $\frac{\partial 
V}{\partial \xi}(\xi, 1) f_1(\xi) = -W_1(\xi)$ and $\frac{\partial V}{\partial \xi}(\xi, 
2) f_2(\xi) = - W_2(\xi)$. It is not hard to see that both pairs $(f_1, W_1)$ and $(f_2, 
W_2)$ have the zero small-time distinguishability property and that the second one is 
not small-time norm-observable. As a matter of fact it is not large-time norm 
observable. In this case Theorem 7 in \cite{Hespanha2} cannot be applied, but according 
to Theorem
\ref{sontagteo}, any $\T^\prime \subseteq \Ta[\tau_D,N_0]$  with fixed, but otherwise 
arbitrary,  $\tau_D>0, N_0\in \N_0$,  is globally uniformly asymptotically stable.

It is worth noting that if one applied Theorem 1 of \cite{Bacciotti}, the attracting set 
so obtained would be the $x_2$-axis (see example 4 of \cite{Bacciotti}). 
          
\section{Invariance for switched systems}
In this section we study the asymptotic behavior of bounded solutions $x$ of the 
switching system (\ref{ss}) corresponding to switching signals which have a positive 
average dwell-time, and  in particular the invariance properties of their $\omega$-limit 
sets. We recall that a point $\xi \in \R^n$ belongs to $\Omega(x)$, the $\omega$-limit 
set of $x:\R_{+} \rightarrow \R^n$, if there exists a strictly increasing sequence of 
times $\{s_k\}$ with $\lim_{k \rightarrow \infty}s_k=\infty$ and $\lim_{k \rightarrow 
\infty}x(s_k)=\xi$. The $\omega$-limit set $\Omega(x)$ is always closed and, when $x$ is 
bounded, it is non-empty, compact and $\lim_{t \rightarrow \infty} d(x(t),\Omega(x))=0$. 
Moreover, $\Omega(x)$ is the smallest closed set which is approached by $x$.

In order to proceed, we will associate to each bounded trajectory $(x,\sigma) \in \Ta$ a 
nonempty subset of $\R^n \times \Gamma$, which we denote $\Omega^{\sharp}(x,\sigma)$, 
and study its invariance properties. 

Let us introduce some more notation and terminology. As stated above, associated with a 
switching signal $\sigma$ there are a strictly increasing sequence of real numbers (the 
sequence of switching times of $\sigma$) $\{t_i\}_{i=0}^{N_\sigma}$, with  $N_{\sigma}$ 
finite or $N_\sigma=\infty$, $t_0=0$ and $\lim_{i \rightarrow \infty}t_i=\infty$ when 
$N_{\sigma}=\infty$, and a sequence of points $\{\gamma_i\}_{i=0}^{N_\sigma} \subseteq 
\Gamma$, with $\gamma_{i}\neq \gamma_{i+1}$ for all $0\le i<N_{\sigma}$, such that 
$\sigma(t) = \gamma_i$ for all $t_i \leq t < t_{i+1}$ with $0 \le i < N_\sigma$, and 
$\sigma(t)=\gamma_{N_{\sigma}}$ for all $t \ge t_{N_\sigma}$ when $N_\sigma$ is finite. 
In order to treat the cases $N_\sigma<\infty$ and $N_\sigma=\infty$ in an unified frame, 
we pick any $\gamma^*\in \Gamma$ and define $t_i=\infty$ and $\gamma_{i}=\gamma^*$ for 
all $i > N_\sigma$ when $N_\sigma$ is finite.

Given a switching signal $\sigma \in \S$ we consider the sequence of maps $\taus^i:\R_+ 
\cup \{\infty\} \rightarrow \R_+ \cup \{\infty\}$, $i \in \N$, defined recursively by:
\begin{itemize} 
\item $\taus^1(t)=t_k$ if $t \in [t_{k-1},t_{k})$ and $\taus^1(\infty)=\infty$;
\item $\taus^{i+1}(t)=\taus^1(\taus^{i}(t))$ for all $t \in \R_+ \cup \{\infty\}$ and 
all $i \ge 2$.
\end{itemize}
We observe that for a given time $t\ge 0$, $\taus^1(t)$ is the first switching time 
greater than $t$, $\taus^2(t)$ is the second switching time greater than $t$, etc. We 
also define, for convenience,  $\taus^0(t)=t$ for all $t \ge 0$.

\begin{defin} \label{omega} Given a bounded trajectory $(x,\sigma) \in \Ta$,
a point $(\xi, \gamma)\in \R^n \times \Gamma$ belongs to $\Omega^\sharp(x,\sigma)$ if 
there exists a strictly increasing and unbounded sequence $\{s_k\} \subset \R_+$ such 
that 
\begin{enumerate}
\item $\lim_{k \rightarrow \infty} \taus^1(s_k)-s_k = r, \; 0<\,r\,\leq \infty$;
\item $\lim_{k\rightarrow \infty}x(s_k)=\xi$ and
$\lim_{k \rightarrow \infty}\sigma(s_k)=\gamma$.\mybox
\end{enumerate} 
\end{defin}

We observe that in the case when $N_\sigma$ is finite, 
$\Omega^\sharp(x,\sigma)=\Omega(x)\times \{\gamma_{N_\sigma}\}$.

The following lemma shows the relation between $\Omega(x)$ and 
$\Omega^{\sharp}(x,\sigma)$. 
\begin{lema} \label{pi1} \rm Let $(x,\sigma)$ be a bounded trajectory belonging to 
$\Ta$. Then $\Omega(x)=\pi_1(\Omega^\sharp(x,\sigma))$.
\end{lema}
{\em Proof.} We only prove the case when $\sigma$ has infinitely many switching times 
since the other case is trivial.

Suppose that $\sigma \in \Sa[\tau_D,N_0]$ for some $\tau_D > 0$ and $N_0 \in \N$ and 
that $\sigma$ has infinitely many switching times. We first note that the inclusion 
$\pi_1(\Omega^\sharp(x,\sigma)) \subseteq \Omega(x)$ readily follows from the definition 
of $\Omega^\sharp(x,\sigma)$. 

In order to prove that $\Omega(x) \subseteq \pi_1(\Omega^\sharp(x,\sigma))$ let $\xi \in 
\Omega(x)$. Then there exists a strictly increasing an unbounded sequence of times 
$\{s_k\}$ such that $x(s_k)\rightarrow \xi$. 

Let $i \ge 0$ be the first integer such that 
\[\limsup_{k \rightarrow \infty} \taus^{i+1}(s_k)-s_k=r, \;{\rm with}\;0<r\le \infty.\]
 Such an integer exists and verifies $i\le  N_0$ since,  due to the definition of 
$\Sa[\tau_D,N_0]$, $\taus^{N_0+1}(s_k)-s_k\ge \tau_D$ for all $k \ge 1$. Let 
$\{s_{k_j}\}$ be a subsequence of $\{s_{k}\}$ so that $\lim_{j\rightarrow \infty} 
\tau^{i+1}_\sigma(s_{k_j})-s_{k_j}=r$. Then, 
i) $\lim_{j\rightarrow \infty} \taus^l(s_{k_j})-s_{k_j}=0$ for all $0 \le l \le i$. 

Consider the sequence $\{\sigma(\taus^{i}(s_{k_j}))\}_{j=1}^{\infty}$; as $\Gamma$ is 
compact, there is a subsequence $\{\sigma(\taus^{i}(s_{k_{j_l}}))\}_{l=1}^{\infty}$ 
which converges to some $\gamma \in \Gamma$. 

We claim that $(\xi,\gamma) \in \Omega^\sharp(x,\sigma)$. In order to prove the claim, 
consider the unbounded sequence  $\{s^\prime_l=\taus^{i}(s_{k_{j_l}})\}$. From i) and 
the facts that $\lim_{j\rightarrow \infty} \taus^{i+1}(s_{k_j})-s_{k_j}=r$ and  
$\taus^1(s^\prime_l)=\taus^1(\taus^{i}(s_{k_{j_l}}))=\taus^{i+1}(s_{k_{j_l}})$, we have 
that $\lim_{l \rightarrow \infty} \taus^1(s^\prime_l)-s^\prime_l=r>0$. 

We note that by construction
 \[\lim_{l\rightarrow \infty}\sigma(s^\prime_l)=\lim_{l\rightarrow 
\infty}\sigma(\taus^{i}(s_{k_{j_l}}))=\gamma.\]

Finally, taking into account that: 
\begin{itemize}
\item $\lim_{l \rightarrow \infty} s^\prime_l-s_{k_{j_l}}=0$;
\item $x$ is uniformly continuous on $\R_{+}$ since $\dot{x}$ is essentially bounded 
($x$ is bounded, $f$ is continuous and $\Gamma$ is compact);
\item $\lim_{k \rightarrow \infty} x(s_k)=\xi$;
\end{itemize}
it follows that $\lim_{l\rightarrow \infty}x(s^\prime_l)=\xi$. Thus $(\xi,\gamma) \in 
\Omega^\sharp(x,\sigma)$ and thereby $\xi \in \pi_1(\Omega^\sharp(x,\sigma))$. Then 
$\Omega(x) \subseteq \pi_1(\Omega^\sharp(x,\sigma))$ and the lemma follows. \square 

The next result shows that, under suitable hypotheses, the set 
$\Omega^{\sharp}(x,\sigma)$ corresponding to a trajectory in $\T_V\cap\T[\chi]\cap 
\Ta[\tau_D,N_0]$ is weakly-invariant w.r.t. $\T_V^*\cap\T[\chi]\cap \Ta[\tau_D,N_0]$. 

\begin{prop} \label{invariance} \rm Let $\chi$ be a closed covering of $\R^n$ and let 
$V:\Omega \times \Gamma \rightarrow \R$ be a continuous function. Suppose that 
$(x,\sigma)$ is a trajectory belonging to $\T_V \cap \T[\chi] \cap \Ta[\tau_D,N_0]$, 
$\tau_D>0, N_0 \in \N$, such that for some compact set $B\subset \Omega$, $x(t)\in B$ 
for all $t\ge 0$. Then $\Omega^\sharp(x,\sigma)$ is weakly-invariant w.r.t. $\T_V^*\cap 
\T[\chi]\cap  \Ta[\tau_D,N_0]$. \mybox
\end{prop}
{\em Proof.} Let $(\xi,\gamma) \in \Omega^\sharp(x,\sigma)$. Then there exists a 
strictly increasing and unbounded sequence $\{s_k\}$ which verifies 1. and 2. of 
Definition \ref{omega}. Let $\sigma_k(\cdot)=\sigma(\cdot+s_k)$ and 
$x_k(\cdot)=x(\cdot+s_k)$. As $\chi$, $V$ and the sequence $\{(x_k,\sigma_k)\}$ are as 
in the hypotheses of Lemma \ref{lemaV} in the Appendix, there exist a subsequence 
$(x_{k_l},\sigma_{k_l})$ and a trajectory $(x^*,\sigma^*)\in 
\T_V\cap\T[\chi]\cap\Ta[\tau_D,N_0]$ such that $\{x_{k_l}\}$ converges uniformly to 
$x^*$ on compact subsets of $\R_{+}$ and $\{\sigma_{k_l}\}$ converges to $\sigma^*$ a.e. 
on $\R^+$. Due to Lemma \ref{cp} in the Appendix, we can also assume without loss of 
generality that $\{\sigma_{k_l}\}$ also verifies condition 2 of that lemma with 
$\sigma^*$ in place of $\sigma$. 

The proof is completed provided we show that $(x^*(0),\sigma^*(0))=(\xi,\gamma)$, 
$(x^*(t),\sigma^*(t))\in \Omega^\sharp(x,\sigma)$ for all $t \ge 0$ and 
$(x^*,\sigma^*)\in \T_V^*$.  

Let us prove first that $(x^*(0),\sigma^*(0))=(\xi,\gamma)$. From the fact that 
$x_k(0)=x(s_k)$, 2. of Definition \ref{omega} and the convergence of $\{x_{k_l}\}$ to 
$x^*$, we have that $x^*(0)=\xi$.

According to 2. of Lemma \ref{cp}, there exists a sequence $\{r_l\} \subset \R_+$ such 
that $\lim_{l \rightarrow \infty}r_l=0$, $\lim_{l \rightarrow 
\infty}\tau^1_{\sigma_{k_l}}(r_l)-r_l>0$ and $\lim_{l \rightarrow 
\infty}\sigma_{k_l}(r_l)=\sigma^*(0)$. From item 1. of Definition \ref{omega} and the 
fact that $\tau^1_{\sigma_k}(0)=\tau_\sigma^1(s_k) - s_k$ for all $k \in \N$, it follows 
that $\lim_{l \rightarrow \infty} \tau^1_{\sigma_{k_l}}(0)>0$. Then $r_l< 
\tau^1_{\sigma_{k_l}}(0)$ for $l$ large enough and, therefore, $\sigma^*(0)=\lim_{l 
\rightarrow \infty}\sigma_{k_l}(r_l)=\lim_{l \rightarrow \infty}\sigma_{k_l}(0)=\lim_{l 
\rightarrow \infty}\sigma(s_{k_l})=\gamma$.

Next we prove that $(x^*(t),\sigma^*(t)) \in \Omega^\sharp(x,\sigma)$ for all $t > 0$. 

Let $t > 0$ and let $\{r_l\}$ be a sequence as in 2. of Lemma \ref{cp}. Consider the 
unbounded sequence $\{s^\prime_l\}$, defined by $s^\prime_l=r_l+s_{k_l}$, which we can 
suppose, without loss of generality, strictly increasing. Due to the fact that 
$\tau^1_{\sigma}(s^\prime_l)=\tau^1_{\sigma_{k_l}}(r_l) + s_{k_l}$, we have that 
$\lim_{l \rightarrow \infty} \tau^1_{\sigma}(s^\prime_l)-s^\prime_l>0$. So 
$\{s^\prime_l\}$ satisfies condition 1. of Definition \ref{omega}.

From 2. of Lemma \ref{cp}, we have that $\lim_{l \rightarrow 
\infty}\sigma(s^\prime_l)=\lim_{l \rightarrow \infty}\sigma_{k_l}(r_l)=\sigma^*(t)$. On 
the other hand, from the uniform convergence of $\{x_{k_l}\}$ to $x^*$ on compact sets 
and the continuity of $x^*$ we have that
$\lim_{l \rightarrow \infty} x(s^\prime_{l})=\lim_{l \rightarrow 
\infty}x_{k_l}(r_l)=x^*(t)$. Hence $(x(s^\prime_{l}),\sigma(s^\prime_{l})) \rightarrow 
(x^*(t),\sigma^*(t))$ as $l \rightarrow \infty$ and thereby $(x^*(t),\sigma^*(t)) \in 
\Omega^\sharp(x,\sigma)$. 

Finally, we prove that $(x^*,\sigma^*)\in \T_V^*$. Since $(x^*,\sigma^*)\in \T_V$ it 
suffices to prove that for any pair of times $t, t^\prime$ with $t < t^\prime$ and 
$\sigma^*(t)=\sigma^*(t^\prime)$, $V(x^*(t),\sigma^*(t))\le 
V(x^*(t^\prime),\sigma^*(t^\prime))$. 

Let $t, t^\prime$ be a pair of times such that $t<t^\prime$ and 
$\sigma^*(t)=\sigma^*(t^\prime)=\gamma$. Since $\lim_{l \rightarrow \infty} 
\sigma_{k_l}(t)=\sigma^*(t)$ a.e. on $\R_+$, and $\sigma^*$ is piecewise constant and 
right continuous, there exists a pair of non-increasing sequences $\{\tau_i\}$, 
$\{\tau^\prime_i\}$ so that 
\begin{itemize}
\item $\tau_i < \tau^\prime_i$ for all $i$;
\item $\tau_i \searrow t$ and $\tau^\prime_i \searrow t^\prime$;
\item for every $i$, $\lim_{l \rightarrow \infty} 
\sigma_{k_l}(\tau_i)=\sigma^*(\tau_i)=\gamma$ and $\lim_{l \rightarrow \infty} 
\sigma_{k_l}(\tau^\prime_i)=\sigma^*(\tau^\prime_i)=\gamma$.
\end{itemize}
Fix $i$. Since $\Gamma$ is finite, there exists $l^*$ such that 
$\sigma_{k_l}(\tau_i)=\sigma_{k_l}(\tau^\prime_i)=\gamma$ for all $l \ge l^*$. Fix 
$l^\prime \ge l^*$. As $\{s_{k_l}\}$ is unbounded, we have that $\tau_i+s_{k_l} > 
\tau_i^\prime+s_{k_{l^\prime}}$ for $l$ large enough, say $l \ge l_0$. Then, for $l \ge 
\max\{l_0,l^*\}$, $\tau_i+s_{k_l} > \tau_i^\prime+s_{k_{l^\prime}}$ and 
$\sigma(\tau_i+s_{k_l})=\sigma_{k_l}(\tau_i)=\sigma_{k_{l^\prime}}(\tau^\prime_i)=\sigma
(\tau^\prime_i +s_{k_{l^\prime}})$. In consequence, 
\begin{eqnarray*}
V(x_{k_{l^\prime}}(\tau_i^\prime), 
\sigma_{k_{l^\prime}}(\tau_i^\prime))=V(x(\tau_i^\prime+s_{k_{l^\prime}}),\sigma(\tau_i^
\prime+s_{k_{l^\prime}})) \ge \\
 V(x(\tau_i+s_{k_{l}}),\sigma(\tau_i+s_{k_{l}}))= V(x_{k_{l}}(\tau_i), 
\sigma_{k_{l}}(\tau_i)).
 \end{eqnarray*}
From the latter, after taking limit as $l \rightarrow \infty$, we get 
$V(x_{k_{l^\prime}}(\tau_i^\prime), \sigma_{k_{l^\prime}}(\tau_i^\prime))\ge 
V(x^*(\tau_i),\sigma^*(\tau_i))$ and from this, letting $l^\prime \rightarrow \infty$, 
we obtain $V(x^*(\tau^\prime_i),\sigma^*(\tau^\prime_i))\ge 
V(x^*(\tau_i),\sigma^*(\tau_i))$. Finally, from the continuity of $V$ and $x^*$ and the 
right continuity of $\sigma^*$, letting $i \rightarrow \infty$ it follows that
$V(x^*(t),\sigma^*(t))\le V(x^*(t^\prime),\sigma^*(t^\prime))$.\square

Now we are ready to prove the extension of LaSalle's invariance principle given in 
Theorem \ref{Lasalle}.

{\bf Proof of Theorem \ref{Lasalle}.}
  It readily follows from Proposition \ref{invariance} and Lemma \ref{omega}. In fact, 
from Proposition \ref{invariance} we have that $\Omega^\sharp(x,\sigma)\subseteq M$. 
Thus, from Lemma \ref{omega}, we deduce that $\Omega(x)\subseteq \pi_1(M)$ and therefore 
that $x(t)$ tends to $\pi_1(M)$ as $t \rightarrow \infty$. \square 
\section{Proof of Theorem \ref{cordin}}
{\em Proof of Theorem \ref{cordin}.}
Since from Theorem \ref{strict} we know that $\T^\prime$ is uniform stable, we only have 
to prove the remaining statements.

Let $(x,\sigma) \in \T^\prime$. Since $\T^\prime$ is uniform stable, it follows that 
$(x,\sigma)$ is bounded and therefore evolves into some compact subset $B$ of 
$\Omega=\R^n$. By applying Theorem \ref{Lasalle} we deduce that $x(t)$ tends to 
$\pi_1(M)=\{0\}$ and the global asymptotic stability of $\T^\prime$ follows.

Suppose now that $\T^\prime \subseteq \Ta[\tau_D,N_0]$. Since $V$ is a weak 
Lyapunov-like function for the family of trajectories 
$\T^*=\T_V\cap\T[\chi]\cap\Ta[\tau_D,N_0]$, and $\T^\prime \subseteq \T^*$, it suffices 
to prove that $\T^*$ is globally uniformly asymptotically stable.

Since we have already proved that $\T^*$ is globally asymptotically stable, and in 
particular uniformly stable, the global uniform asymptotic stability of $\T^*$ will be 
established if we show that $\T^*$ verifies (a) of Remark \ref{gas}. 

As $\T^*$ is invariant by time translations, i.e., for all $s \ge 0$, 
$(x(\cdot+s),\sigma(\cdot+s)) \in \T^*$ if $(x,\sigma) \in \T^*$, in order to prove (a) 
of Remark \ref{gas} it is sufficient to  
show that $\T^*$ verifies the weaker condition:
\begin{itemize}
\item[(*)] for each $R>0$ and each $\varepsilon > 0$, there exists
$T\ge 0$ such that for all $(x, \sigma) \in \T^*$, 
\[ |x(0)| < R \Longrightarrow |x(t)|< \varepsilon \quad \forall t \ge T.\]
\end{itemize}

Suppose that (*) does not hold. Then there exist $\varepsilon_0>0$, $\eta_0>0$, a 
sequence of trajectories $\{(x_k,\sigma_k)\} \subset \T^*$ and an increasing and 
unbounded sequence of times $\{\tau_k\}$ such that $|x_k(0)|\le \eta_0$ and 
$|x_k(\tau_k)|\ge \varepsilon_0$ for all $k$. 

Since $\{(x_k,\sigma_k)\}$ is uniformly bounded, from Lemma \ref{lemaV} we know that 
there exists a subsequence $\{(x_{k_l},\sigma_{k_l})\}$ and a trajectory $(x^*,\sigma^*) 
\in \T^*$ such that $\{x_{k_l}\}$ converges to $x^*$ uniformly on compact sets. Let 
$\varepsilon^\prime=\alpha^{-1}(\varepsilon_0/2)$, with $\alpha$ as in 1. of Definition 
\ref{stability}. As $x^*$ converges to $0$ as $t \rightarrow \infty$, there exists a 
time $T >0$ such that $|x^*(T)|< \varepsilon^\prime$. Since $\{x_{k_l}\}$ converges to 
$x^*$ uniformly on compact sets, $|x_{k_l}(T)| < \varepsilon^\prime$ for $l$ large 
enough. Then, due to the uniform stability of $\T^*$, we have that, for $l$ large enough 
and $t\ge T$, 
\[ |x_{k_l}(t)| \le \alpha(|x_{k_l}(T)|) \le 
\alpha(\varepsilon^\prime)=\frac{\varepsilon_0}{2},\]  
which is a contradiction. \square
\section{Conclusions}
In this paper we have presented an extension of LaSalle's invariance principle for 
switched nonlinear systems assuming that the family of subsystems is finite and that the 
switching signals have a positive average dwell-time. This extension enabled us to 
obtain some asymptotic stability criteria for this class of systems. Examples were 
presented that show the application  of our results to cases either intractable with 
some of the previously mentioned results or upon which those results give no conclusive 
answers. In addition, results about the compactness of the trajectories of the systems 
involved were exhibit that not only were instrumental used in the proof of some results 
along the paper, but are important by their own.

Finally, we point out that extensions of LaSalle's principle for switched nonlinear 
systems in the case that $\Gamma$ is infinite and also integral invariance principles 
for the same class of switched systems have been already obtained and are currently 
under preparation for publication.

\bigskip

\appendix

\noindent
{\Large{\bf Appendix}}

\section{Some compactness results of trajectories of switched systems}
In this Appendix we will show that, under suitable hypotheses, certain families of 
trajectories of system (\ref{ss}) enjoy a certain kind of sequential compactness. 
 
We say that a sequence $\{(x_k,\sigma_k)\}$ of trajectories of (\ref{ss}) is uniformly 
bounded if there exists $M\ge 0$ such that for all $k$, $|x_k(t)| \le M$ for all $t \ge 
0$.  
\begin{defin} \rm  A family $\T^*$ of maximal trajectories of (\ref{ss}) has the {\bf 
SC} property if for every uniformly bounded sequence $\{(x_k,\sigma_k)\}\subset \T^*$ 
there exist a subsequence $\{(x_{k_l},\sigma_{k_l})\}$ and a trajectory 
$(x^*,\sigma^*)\in \T^*$ such that $\{x_{k_l}\}$ converges to $x^*$ uniformly on compact 
sets of $\R_{+}$ and $\lim_{l \rightarrow \infty}\sigma_{k_l}(t)=\sigma^*(t)$ a.e. on 
$\R_{+}$. \mybox
\end{defin}
\begin{prop} \rm \label{SC} Assume that $\chi$ is a closed covering of $\R^n$.Then 
$\T[\chi] \cap \Ta[\tau_D,N_0]$ has the {\bf SC} property for all $\tau_D>0$ and  all 
$N_0 \in \N$. \mybox
\end{prop}

The following lemma is used in the proof of Proposition \ref{SC} and in some parts of 
section 4. 
\begin{lema} \rm \label{cp} Let $\{\sigma_k\}$ be a sequence of switching signals in 
$\Sa[\tau_D,N_0]$ with $\tau_D>0$ and $N_0 \in \N$.

 Then there exist a subsequence $\{\sigma_{k_l}\}$ and a switching signal $\sigma \in 
\Sa[\tau_D,N_0]$ such that 
\begin{enumerate}
\item 
$\lim_{l \rightarrow \infty}\sigma_{k_l}(t)=\sigma(t)$ for almost all $t\ge 0$;
\item for each $t \in \R_{+}$ there exists a sequence of positive times 
$\{r_l\}_{l=1}^\infty$ such that $\lim_{l\rightarrow \infty}r_l=t$,  
\[ \lim_{l\rightarrow \infty}\sigma_{k_l}(r_l)=\sigma(t)\quad {\rm and}\quad 
\lim_{l\rightarrow \infty} \tau^1_{\sigma_{k_l}}(r_l)-r_l>0.\]
\end{enumerate} \mybox
\end{lema}
{\em Proof.}  Let $\Rc= \R_{+} \cup \{\infty\}$ be the one-point compactification of 
$\R_{+}$, which we recall is a compact metric space, and let $\mathcal{K}=(\Rc \times 
\Gamma)^{\N_0}$ be the set of all the sequences $p=\{(t_i,\gamma_i):\; t_i\in \Rc,\; 
\gamma_i \in \Gamma, i \in \N_0 \}$ endowed with the product topology. We note that 
since $\mathcal{K}$ is the Cartesian product of a countable number of compact metric 
spaces, it is metrizable (see \cite{Dugundji}, Theorem 7.2 on p. 190) and compact (see 
\cite{Dugundji}, Theorem 1.4 on p. 224 ).

For each $k \in \N$, let $\{t_i^k\}_{i=0}^\infty$ be the sequence of switching times 
associated to the switching signal $\sigma_k$ and let $\{\gamma_i^k\}_{i=0}^\infty$ be 
the sequence of points of $\Gamma$ defined by $\gamma_i^k=\sigma_k(t_i^k)$. Observe that 
when $N_{\sigma_k}$ is finite we have, according to the convention above, that 
$\gamma_i^k = \gamma^*$ and $t_i^k = \infty$ for every $i > N_{\sigma_k}$. Let 
$p_k=\{p_k(i)=(t^k_i,\gamma^k_i)\}_{i=0}^{\infty} \in \mathcal{K}$.

As $\mathcal{K}$ is a compact metric space, there exists a subsequence 
$\{p_{k_l}\}_{l=1}^{\infty}$ which converges, say to 
$p=\{p(i)=(t_i,\gamma_i)\}_{i=0}^{\infty}$, i.e., for each $i \in \N_0$, $t^{k_l}_i 
\rightarrow t_i$ and $\gamma^{k_l}_i \rightarrow \gamma_i$ as $l \rightarrow \infty$. 

Since for each $k \in \N$, $\{t^k_i\}_{i=0}^{\infty}$ is nondecreasing and $t^k_0=0$, it 
readily follows that $\{t_i\}_{i=0}^{\infty}$ is nondecreasing and $t_0=0$.

We claim that:
\begin{itemize}
\item[(a)] For every open interval $(a,b)$, with $a<b$, the number of indexes $i\in \N$ 
such that $t_i \in (a,b)$, is bounded by $N_0 +(b-a)/\tau_D$;
\item[(b)] the number of indexes $i$ such that $t_i=0$ is at most $N_0+1$. 
\end{itemize}
We will only prove (a) since (b) can be proved in a similar way.\\
{\em Proof of (a).} Suppose on the contrary that there are $r>N_0+(b-a)/\tau_D$ indexes, 
say $i_1,\ldots,i_r$, such that $t_{i_j}\in (a,b)$ for $j = 1,\ldots, r$. Since  
$t^{k_l}_{i_j} \rightarrow t_{i_j}$ for each $j=1,\ldots,r$, we have that $t^{k_l}_{i_j} 
\in (a,b)$ for all $j=1,\ldots,r$ if $l$ is large enough, which contradicts the fact 
that $\sigma_{k}$ belongs to $\Sa[\tau_D,N_0]$.

In order to define $\sigma$ as in the thesis of the lemma, let $\{i_j\}_{j=0}^{N}$, with 
$N \le \infty$, be the unique subsequence of $\N_0$ that verifies: 
\begin{itemize}
\item $0=t_0=\cdots=t_{i_0}$ and $t_{i_0+1}>0$;
\item $t_{i_{j}+1}=\cdots=t_{i_{j+1}}$ and $t_{i_{j+1}}<t_{i_{j+1}+1}$ for all
$0 \le j < N$;
\item $t_{i_N+1}=\infty$ when $N<\infty$. 
\end{itemize}
Note that this subsequence is well defined due to (a) and (b). 

Now, let $\sigma:\R_{+}\rightarrow \Gamma$ be the switching signal defined by: 
$\sigma(t)=\gamma_{i_j}$ for all $t \in [t_{i_j},t_{i_{j+1}})$ and all $j \le N$. From 
(a) it readily follows that $\sigma \in \Sa[\tau_D,N_0]$. 

Now we proceed to prove 1. and 2. of the thesis of the lemma. We will consider two 
cases.

{\em Case I.} $t \notin \{t_{i_j}\}_{j=0}^N$. 

Let $j^* \in \N_0$ so that $t_{i_{j^*}}<t < t_{i_{j^*}+1}$. As $\lim_{l \rightarrow 
\infty}t^{k_l}_{i_{j^*}}=t_{i_{j^*}}$ and $\lim_{l \rightarrow \infty 
}t^{k_l}_{i_{j^*}+1} = t_{i_{j^*}+1}$, $t^{k_l}_{i_{j^*}}<t<t^{k_l}_{i_{j^*}+1}$ for $l$ 
large enough, say $l \ge L$. Therefore, for $l \ge L$, 
$\sigma_{k_l}(t)=\sigma_{k_l}(t^{k_l}_{i_{j^*}})=\gamma^{k_l}_{i_{j^*}}$ and, 
consequently, 
\begin{eqnarray}\label{I}
\lim_{l \rightarrow \infty}\sigma_{k_l}(t)=\lim_{l \rightarrow 
\infty}\gamma^{k_l}_{i_{j^*}}=\gamma_{i_{j^*}}=\sigma(t).
\end{eqnarray} 
As $\{t_{i_j}\}_{j=0}^N$ is a set of measure zero, the latter shows 1).

From the fact that $t^{k_l}_{i_{j^*}}<t<t^{k_l}_{i_{j^*}+1}$ for $l \ge L$, we also have 
that
\begin{eqnarray*}
\lim_{l \rightarrow \infty} \tau^1_{\sigma_{k_l}}(t)-t=\lim_{l \rightarrow \infty} 
t^{k_l}_{i_{j^*}+1}-t= t_{i_{j^*}+1}-t>0,   
\end{eqnarray*}
which shows that 2) holds with the sequence $\{r_l\}$ defined by $r_l=t$ for all $l \in 
\N$. 

{\em Case II.} $t=t_{i_{j^*}}$ for some $j^* \in \N_0$. 

By using arguments similar to those used in the preceding case, we have that 
 \begin{equation*}
\lim_{l \rightarrow \infty} \tau^1_{\sigma_{k_l}}(t)-t^{k_l}_{i_{j^*}}=\lim_{l 
\rightarrow \infty} t^{k_l}_{i_{j^*}+1}-t^{k_l}_{i_{j^*}}= t_{i_{j^*}+1}-t_{i_{j^*}}>0,   
\end{equation*}
and that $\lim_{l \rightarrow \infty} \sigma(t^{k_l}_{i_{j^*}})=\sigma(t)$.

Consequently, item 2. holds with the sequence $\{r_l\}$ defined by 
$r_l=t^{k_l}_{i_{j^*}}$ for all $l \in \N$. 
\square \\
\begin{remark}\label{comgamma}
It is worth mentioning that it is not necessary that $\Gamma$ be a finite set for the 
thesis of Lemma \ref{cp} to hold. In fact, as can be easily seen from its proof, it 
suffices that $\Gamma$ be a compact metric space. \mybox
\end{remark}   
\begin{remark}\label{kloeden}
A result on compactness of switching signals,  proved with arguments different to ours, 
has recently appeared in \cite{Kloeden}. That result (Theorem 1 of that paper) states 
that given $\Gamma = \{1, \dots, N\}$, the set of switching signals with a {\em fixed} 
dwell-time $\tau_D > 0$  is a compact subset of the metric space $(\S, d)$ where the 
metric $d$ is defined by 
\begin{eqnarray*}
d(u, v) = \sum_{n = 1}^\infty 2^{-n} \int_0^n |u(s) - v(s)| ds, 
\end{eqnarray*}
for all $u$ and $v$ in $\S$.  

By using  Lemma \ref{cp} one can generalize  Theorem 1 in \cite{Kloeden} to average 
dwell-time signals (the {\em fixed} dwell-time hypothesis is essential in the proof 
given in \cite{Kloeden}). In fact, if we consider the metric space $(\Sa[\tau_D, N_0], 
d)$ with $d$ the metric above, the compactness of $\Sa[\tau_D, N_0]$ follows from Lemma 
\ref{cp} and the application of Lebesgue's Dominated Convergence Theorem. Moreover, 
since as 
pointed out in Remark \ref{comgamma}, Lemma \ref{cp} holds for a compact metric space 
$(\Gamma, \rho)$, $\Sa[\tau_D, N_0]$ is compact with the metric
\begin{eqnarray*}
d(u, v) = \sum_{n = 1}^\infty 2^{-n} \int_0^n \rho(u(s), v(s)) ds, 
\end{eqnarray*}
for all $u$ and $v$ in $\S$. \mybox  
\end{remark}   
{\em Proof of Proposition \ref{SC}.} As $\{(x_k,\sigma_k)\}$ is bounded there exists $M 
\ge 0$ such that $|x_k(t)|\le M$ for all $t\ge 0$ and all $k \in \N$. Let 
$M^\prime=\max_{|\xi| \le M, \gamma \in \Gamma }|f(\xi,\gamma)|$. Then, for every 
positive integer $k$, $|\dot{x}_k(t)| \le M^\prime$ for almost all $t \in [0,+\infty)$. 
In consequence $\{x_k\}$ is equibounded and equicontinuous. Then, applying the 
Arzela-Ascoli Theorem we deduce the existence of a subsequence $\{x_{k_l}\}$ and a 
continuous function $x^*:\R_{+} \rightarrow \R^n$ such that $\{x_{k_l}\}$ converges to 
$x^*$ uniformly on compact subsets of $\R_{+}$.

Consider the subsequence $\{\sigma_{k_l}\}$. Due to Lemma \ref{cp} we can suppose 
without loss of generality that there exists $\sigma^* \in \Sa[\tau_D,N_0]$ such that 
$\lim_{l \rightarrow +\infty}\sigma_{k_l}(t)=\sigma^*(t)$ for almost all $t\ge 0$.

We claim that $x^*$ is a solution of (\ref{ss}) corresponding to $\sigma^*$ and, in 
consequence, $(x^*,\sigma^*)\in \Ta[\tau_D,N_0]$.

Let $t>0$. Then 
\begin{eqnarray*}
x^*(t) = \lim_{l\rightarrow +\infty} x_{k_l}(t) &=& 
\lim_{l \rightarrow +\infty}\left (x_{k_l}(0)+\int_0^t f(x_{k_l}(s),\sigma_{k_l}(s)) ds 
\right )\\
&=& x^*(0)+ \lim_{l \rightarrow +\infty} \int_0^t f(x_{k_l}(s),\sigma_{k_l}(s)) ds.
\end{eqnarray*}
As $|f(x_{k_l}(s),\sigma_{k_l}(s))|\le M^\prime$ for all $s \in [0,t]$ and $\lim_{l 
\rightarrow +\infty}f(x_{k_l}(s),\sigma_{k_l}(s))=f(x^*(s), \sigma^*(s))$ for almost all 
$s \in [0,t]$, applying the Lebesgue Dominated Convergence Theorem we have that 
\[  \lim_{l \rightarrow +\infty} \int_0^t f(x_{k_l}(s),\sigma_{k_l}(s)) ds= 
\int_0^t f(x^*(s),\sigma^*(s)) ds, \]
and, consequently,
\[ x^*(t)=x^*(0)+\int_0^t f(x^*(s),\sigma^*(s))ds.\]

It only remains to show that $(x^*,\sigma^*)$ belongs to $\T[\chi]$.
Let $t \ge 0$ so that $\lim_{l\rightarrow \infty}\sigma_{k_l}(t)=\sigma^*(t)=j$. As 
$\Gamma$ is finite, there exists $l^*>0$ such that $\sigma_{k_l}(t)=j$ for all $l \ge 
l^*$. Note that due to the definition of $\T[\chi]$, we also have that $x_{k_l}(t) \in 
\chi_j$ for all $l \ge l^*$. As $\{x_{k_l}(t)\}$ converges to $x^*(t)$ and $\chi_j$ is 
closed, we deduce that $x^*(t) \in \chi_j$. In consequence, 
\begin{eqnarray} \label{a.a.}
 x^*(t) \in \chi_{\sigma^*(t)}\quad {\rm for \; almost \;all}\; t \ge 0.
\end{eqnarray}
Now let $t\ge 0$ be arbitrary. Due to (\ref{a.a.}) there exists a sequence  $\{s_k\}$ 
which converges to $t$ so that $t \le s_k$, $\sigma^*(s_k)=\sigma^*(t)$ and $x^*(s_k)\in 
\chi_{\sigma^*(t)}$ for all $k$. Then, from the continuity of $x^*$ and the fact that 
$\chi_{\sigma^*(t)}$ is closed, we have that $\lim_{k\rightarrow 
\infty}x^*(s_k)=x^*(t)\in \chi_{\sigma^*(t)}$ and the proof is completed. \square
\begin{lema}\rm \label{lemaV} Assume that $\chi$ is a closed covering of $\R^n$ and  
$V:\Omega \times \Gamma \rightarrow \R$, with $\Omega$ an open subset of $\R^n$, is 
continuous. Let $\{(x_k,\sigma_k)\}$ be a sequence of maximal trajectories of (\ref{ss}) 
belonging to $\T_V \cap \T[\chi] \cap \Ta[\tau_D,N_0]$, $\tau_D>0$, $N_0 \in \N$, and 
suppose that there exists a compact subset $B \subset \Omega$ such that $x_k(t) \in B$ 
for all $t \ge 0$ and all $k$.

Then there exist a subsequence $\{(x_{k_l},\sigma_{k_l})\}$ and a maximal trajectory 
$(x^*,\sigma^*)\in \T_V \cap \T[\chi] \cap \Ta[\tau_D,N_0]$ such that $\{x_{k_l}\}$ 
converges to $x^*$ uniformly on compact sets of $\R_+$ and $\lim_{l \rightarrow \infty} 
\sigma_{k_l}(t)=\sigma^*(t)$ a.e. on $\R_+$.\mybox
\end{lema}     
{\em Proof.} Since $\{(x_{k},\sigma_{k})\}$ is uniformly bounded, from Proposition 
\ref{SC} there exist a subsequence $\{(x_{k_l},\sigma_{k_l})\}$ and a maximal trajectory 
$(x^*,\sigma^*)\in \T[\chi] \cap \Ta[\tau_D,N_0]$ such that $\{x_{k_l}\}$ converges to 
$x^*$ uniformly on compact sets of $\R_+$ and $\lim_{l \rightarrow \infty} 
\sigma_{k_l}(t)=\sigma^*(t)$ a.e. on $\R_+$. Therefore the lemma follows provided 
$(x^*,\sigma^*)\in \T_V$. As $x^*(t) \in B \subset \Omega$ for all $t \in \R_{+}$, it 
only remains to prove that for all pair of times $t, t^\prime$, with $t<t^\prime$ and 
$\sigma^*(t)=\sigma^*(t^\prime)$, $V(x^*(t),\sigma^*(t))\ge 
V(x^*(t^\prime),\sigma^*(t^\prime))$. 

Let $t, t^\prime$ be a pair of times such that $t<t^\prime$ and 
$\sigma^*(t)=\sigma^*(t^\prime)=\gamma$. Since $\lim_{l \rightarrow \infty} 
\sigma_{k_l}(t)=\sigma^*(t)$ a.e. on $\R_+$, and $\sigma^*$ is piecewise constant and 
continuous from the right, there exists a pair of non-increasing sequences $\{s_i\}$, 
$\{s^\prime_i\}$ so that 
\begin{itemize}
\item $s_i < s^\prime_i$ for all $i$;
\item $s_i \searrow t$ and $s^\prime_i \searrow t^\prime$;
\item for every $i$, $\lim_{l \rightarrow \infty} 
\sigma_{k_l}(s_i)=\sigma^*(s_i)=\gamma$ and  $\lim_{l \rightarrow \infty} 
\sigma_{k_l}(s^\prime_i)=\sigma^*(s^\prime_i)=\gamma$.
\end{itemize}
Fix $i$. Since $\Gamma$ is finite,  we have that 
$\sigma^*_{k_l}(s_i)=\sigma^*_{k_l}(s^\prime_i)=\gamma$ for $l$ large enough, say $l \ge 
l^*$. Then, from the definition of $\T_V$ and the fact that $s_i < s^\prime_i$, it 
follows that $V(x_{k_l}(s_i),\sigma_{k_l}(s_i))\ge 
V(x_{k_l}(s^\prime_i),\sigma_{k_l}(s^\prime_i))$ for all $l \ge l^*$. In consequence, by 
taking limit as $l \rightarrow \infty$ we obtain that $V(x^*(s_i),\sigma^*(s_i))\ge 
V(x^*(s^\prime_i),\sigma^*(s^\prime_i))$ and, {\em a posteriori}, letting $i \rightarrow 
\infty$ we have that $V(x^*(t),\sigma^*(t))\ge 
V(x^*(t^\prime),\sigma^*(t^\prime))$.\square


\begin{thebibliography}{99}

\bibitem{Bacciotti} A. Bacciotti and L. Mazzi, An invariance principle for nonlinear 
switched systems,   {\em Rapporti interni Nr. 22}, Departimento di Matematica, 
Politecnico di Torino, 2004.  

\bibitem{Branicky} M. Branicky, Multiple Lyapunov functions and other analysis tools for 
switched and hybrid systems, {\em IEEE Trans. Automat. Contr.}, vol.  43, 475-482, 1998.

\bibitem{DeCarlo} R.A. DeCarlo, M.S. Branicky, S. Pettersson and B. Lennartson,  
Perspectives and results on the stability and stabilizability of hybrid systems, {\em 
Proc. IEEE}, vol. 88, pp. 1069-1082, 2000.

\bibitem{Dugundji} J. Dugundji, {\em Topology}, Allyn and Bacon, Inc., Boston, 1966. 

\bibitem{Hespanha} J.P Hespanha, Uniform stability properties of switched linear 
systems: extensions of LaSalle's Invariance Principle, {\em IEEE Trans. Automat. 
Control}, vol. 49, pp. 470-482, 2004.

\bibitem{Hespanha2} J.P. Hespanha, D. Liberzon, D. Angeli and E. Sontag, Nonlinear 
observability notions and stability of switched systems, {\em IEEE Trans. Automat. 
Control}, vol. 50, pp. 154-168, 2005.

\bibitem{Kloeden} P.E. Kloeden, Nonautonomous attractors of switched systems, submitted, 
2004.
\bibitem{lasalle} J. P. LaSalle, {\em The Stability of Dynamical Systems}, ser. Regional 
Conference Series in Applied Mathematics, Philadelphia, PA: SIAM, 1996.
\bibitem{Liberzonbook} 
D. Liberzon, {\em Switching in Systems and Control},  
Birkh\"{a}user, Boston, 2003.

\bibitem{liberzon} D. Liberzon and A.S. Morse, Basic problems in stability and design of 
switched systems, {\em IEEE Control Systems Magazine}, vol. 19, pp. 59-70, 1999.

\bibitem{Lin} Y. Lin, E.D. Sontag and Y. Wang, A smooth converse Lyapunov theorem for 
robust stability, {\em SIAM J. Control Optim.}, vol. 34, pp.  124-160, 1996.

\bibitem{lygeros} J. Lygeros, K. H. Johansson, S. M. Simi\'{c}, J. Zhang and S. Sastry, 
Dynamical Properties of Hybrid Automata, {\em IEEE Trans. Automat. Contr.} vol. 48, pp. 
2-17, 2003.

\bibitem{Peleties} P. Peleties and R.A. DeCarlo, Asymptotic stability of $m$-switched 
systems using Lyapunov-like functions, in {\em Proc. 1991 Amer. Control Conf.}, 1991. 
pp. 1679-1684.

\end{thebibliography}
\end{document}